\documentclass[12pt]{amsart}
\theoremstyle{plain}
\newtheorem{thm}{Theorem}[section]
\newtheorem{theorem}[thm]{Theorem}

\newtheorem{proposition}[thm]{Proposition}
\theoremstyle{definition}

\newtheorem{definition}[thm]{Definition}

\newtheorem{problem}[thm]{Problem}

\newtheorem{question}[thm]{Question}

\numberwithin{equation}{section}

\newcommand{\sC}{{\mathcal C}}
\newcommand{\sD}{{\mathcal D}}

\newcommand{\sF}{{\mathcal F}}

\newcommand{\sK}{{\mathcal K}}

\newcommand{\sO}{{\mathcal O}}
\newcommand{\sP}{{\mathcal P}}

\newcommand{\sT}{{\mathcal T}}

\newcommand{\sW}{{\mathcal W}}

\newcommand{\C}{{\mathbb C}}

\newcommand{\BP}{{\mathbb P}}

\begin{document}

\title[VMRT-structures]{Mori geometry meets Cartan geometry: \\ Varieties of minimal rational tangents}
\author[Jun-Muk Hwang]{Jun-Muk Hwang${}^1$}
\address{Korea Institute for Advanced Study, Hoegiro 87, Seoul, 130-722, Korea}
\email{jmhwang@kias.re.kr}
\thanks{${}^1$ supported
by National Researcher Program 2010-0020413 of NRF}

\begin{abstract}
We give an introduction to the theory of varieties of minimal rational tangents, emphasizing its aspect as a fusion of algebraic geometry and differential geometry, more specifically, a fusion of Mori geometry of minimal rational curves and Cartan geometry of cone structures.
\end{abstract}

\maketitle
\section{Introduction: a brief prehistory}
Lines have been champion figures in classical geometry.
Together with circles, they  dominate the entire geometric contents of  Euclid. Their dominance is no less strong in projective
geometry. Classical projective geometry is full of fascinating results about intricate combinations of lines.    As geometry entered the modern era, lines  evolved into objects of greater flexibility and generality while retaining
all the beauty and brilliance of classical lines. As  Euclidean geometry developed  into Riemannian  geometry, for example,  lines were replaced by  geodesics which then inherited all the glory of
Euclidean lines.

In the transition from   classical projective geometry to complex projective geometry, real lines have been replaced by complex lines.     Lines over complex numbers  have all  the power of lines in
classical projective geometry and even more: results of greater elegance and harmony are obtained over complex numbers.
A large number of results on lines and their interactions with other varieties have been obtained in complex projective geometry, their dazzling beauty no less impressive than that of classical geometry.
 But as complex projective geometry develops further into  complex geometry and abstract algebraic geometry,
which emphasize intrinsic properties of complex manifolds and abstract varieties, the notion of lines in projective space
seems to be too limited for it to keep its leading role.

Firstly, to be useful in intrinsic geometry of projective varieties in projective space, lines should lie on the projective varieties. But most projective varieties do not contain lines. Even when a projective manifold contains  lines, the locus of lines is, often, small and then such a locus is usually regarded as an exceptional part.
Of course,  there are many important varieties that are covered by lines, but they belong to a limited class from the general perspective of classification theory of  varieties. In short,
\begin{quote}
(*)\; \em the class of projective manifolds covered by lines seems to be too special from the perspective of the general theory of complex
manifolds or algebraic varieties.
\end{quote}
  Secondly, many of the methods employed to use lines on varieties in projective space depend on the extrinsic geometry
   of ambient projective space. They do not truly belong to intrinsic geometry of the varieties.   Such geometric arguments  are undoubtedly useful in fathoming deeper geometric properties of varieties which are  described explicitly,  at least to some extent.  But can such methods
   yield  results on {\em a priori} unknown varieties, defined abstractly by intrinsic conditions?      In short,
\begin{quote}
  (**) \; \em tools employed in  line geometry are not intrinsic enough to handle intrinsic problems on abstractly described varieties.
\end{quote}
 These concerns show that lines in projective space have  a rather limited role in the modern development of complex algebraic geometry. Is there a  more general and more powerful notion in complex algebraic geometry that can replace the role of lines, as
   geodesics do  in Riemannian geometry? No serious candidate had emerged  until Mori's groundbreaking work \cite{Mori}.

 In the celebrated paper \cite{Mori}, Mori shows that
 a large class of projective manifolds, including all Fano manifolds, are covered by certain intrinsically defined rational curves  that behave like lines in many respects.  Let us call these rational curves `minimal rational curves'. If a projective manifold embedded in projective space is covered by lines,
 these lines are  minimal rational curves of the projective manifold,  so the notion of minimal rational curves can be viewed as an intrinsic generalization of lines.

 The class of projective manifolds covered by  minimal rational curves are called uniruled projective manifolds.  Generalizing Mori's result, Miyaoka and Mori  have proved in \cite{MiyaokaMori} that a projective manifold is uniruled if its anti-canonical bundle satisfies a certain positivity condition. This implies that uniruled projective manifolds form a large class of algebraic varieties.  Furthermore, the minimal model program, a modern structure theory of higher-dimensional
 algebraic varieties, predicts that uniruled projective manifolds are precisely those projective manifolds that do not admit minimal models. Thus \begin{quote} \em projective manifolds covered by minimal rational curves form a distinguished class of manifolds, worthy of independent study from the view-point of classification theory of  general projective varieties, and at the same time,  large enough
   to contain examples of great diversity. \end{quote}  This overcomes the limitation (*) of the class of projective manifolds covered by lines.

 Furthermore, Mori's work exhibits how to use minimal rational curves in an intrinsic way to obtain geometric information on uniruled projective manifolds. The main tool here is the deformation theory of curves, a machinery of modern complex algebraic geometry- somewhat reminiscent of the use of variational calculus in the local study of geodesics in Riemannian geometry. An example is the property that  a minimal rational curve cannot be deformed when two distinct points on the curve are  fixed. This result  generalizes   the fundamental postulate of classical geometry that "two points determine one line".  The important point is that such a classical property of lines can be recovered by modern deformation theory in an abstract setting.  \begin{quote} \em  Deformation theory of rational curves is a powerful technique applicable to conceptual  problems on varieties defined in abstract intrinsic terms. \end{quote}    In \cite{Mori}, Mori, in fact,  has resolved one of the toughest problems of this kind,  the Hartshorne conjecture,  characterizing projective space by the positivity of the tangent bundle. The methods employed in the theory of minimal rational curves are certainly free of the concern (**) on the
 tools of line geometry.

These considerations indicate that minimal rational curves can serve as the natural generalization of lines,  overcoming the limitations of lines, while inheriting  their powerful and elegant features.

\medskip
Our main interest is the geometry of minimal rational curves in uniruled projective manifolds.  As in Mori's work, we would like to see how minimal rational curves can be used to control the intrinsic geometry of uniruled projective manifolds.   One guiding problem is  the following question on recognizing a given uniruled projective manifolds by  minimal rational curves.

 \begin{problem}\label{p.recognition}  Let $S$ be  a (well-known) uniruled projective manifold. Given another uniruled projective manifold $X$, what properties of
 minimal rational curves on $X$ guarantee that $X$ is biregular (i.e. isomorphic as abstract algebraic varieties)  to  $S$? \end{problem}

Here, the setting of the problem is algebraic geometry and by properties of minimal rational curves, we mean algebro-geometric
properties.
 When $S$ is projective space, a version of this problem is precisely what Mori solved in \cite{Mori}.
 The initial goal of \cite{Mori} was to   prove the Hartshorne conjecture, which characterizes projective space by certain positivity
 property of the tangent bundle. After showing that the projective manifold in question is uniruled, Mori used the tangent directions of minimal rational curves to finish the proof.  This part of Mori's proof has been greatly strengthened by the later work of Cho-Miyaoka-Shepherd-Barron \cite{CMSB}, which says roughly  the following (see Theorem \ref{t.CMSB2} for a precise statement).
\begin{theorem}\label{t.CMSB}  Suppose for a general point $x$ on a uniruled projective manifold $X$ and a general tangent direction
  $\alpha \in \BP T_x(X)$, there exists a minimal rational curve through $x$ tangent to $\alpha$. Then
  $X$ is projective space. \end{theorem}
This is a very satisfactory answer to Problem \ref{p.recognition}  when $S$ is projective space. It includes, as special cases, many
previously known characterizations of projective space.  One may wonder why the condition on minimal rational curves here is formulated in terms of their tangent directions, not in terms of some other properties of minimal rational curves. The essential reason is because  the main technical tool to handle minimal rational curves is the deformation theory of curves, as mentioned before.  The
  tangential information of curves is essential in deformation theory.
For this reason, it is natural and also useful to give conditions in terms of tangent directions of minimal rational curves.

  What about other uniruled projective manifolds? When $S$ is different from projective space, Theorem \ref{t.CMSB} says
  that   minimal rational curves on $S$ exist only in some distinguished directions. Thus in the setting of Problem \ref{p.recognition},  it is natural to consider
\begin{quote} \em
  the subvariety  $\sC_s \subset \BP T_s(S)$ consisting of the directions of minimal rational curves
 through $s \in S$ and the corresponding subvariety $\sC_x \subset \BP T_x(X)$.
\end{quote}
 When $S$ is projective space, we have $\sC_s = \BP T_s(S)$ for any $s \in S$. Theorem \ref{t.CMSB} says that a uniruled
 projective manifold $X$ is projective space if and only if $\sC_x = \BP T_x(X)$ at some point $x \in X$.  In other words, if a uniruled projective manifold has the same type of $\sC_x$ as projective space, then it is projective space. Based on this observation,    we can
 refine our guiding Problem \ref{p.recognition} as follows.
 \begin{problem}\label{p.VMRT}  Let $S$ be  a (well-known) uniruled projective manifold. Given another uniruled projective manifold $X$, what properties of $\sC_x \subset \BP T_x(X)$ for general points $x \in X$
 guarantee that $X$ is biregular   to  $S$? \end{problem}

Comparing this with Theorem \ref{t.CMSB}, one may wonder why we are asking for information on $\sC_x$ for general points $x \in X$,
instead of a single point $x \in X$ as in Theorem \ref{t.CMSB}. This is because the information at one point $x \in X$ seems to be too weak to characterize $X$ when  $\sC_x \neq \BP T_x(X).$ The equality $\sC_x = \BP T_x(X)$ implies that minimal rational curves through one point $x$ cover the whole of $X$. This is why in Theorem \ref{t.CMSB} the information at one point is sufficient to control the whole of $X$. If $\sC_x \neq \BP T_x(X)$,  minimal rational curves through one point $x$ cover only small part of $X$.
Besides,  the subvariety $\sC_s \subset \BP T_s(S)$
 may  change  as the  point $s \in S$ varies and so  the expected condition is not just on  $\sC_x$ for a single $x$, but   on the family $$\{ \sC_x \subset \BP T_x(X), \mbox{ general } x \in X \}.$$
This is why we are asking for the data $\sC_x$  for all general $x \in X$.

Now as in Problem \ref{p.recognition}, the properties of $\sC_x$ that we are looking for in Problem \ref{p.VMRT}  are algebro-geometric properties.
   In algebraic geometry, however, to use properties of  such a family of varieties to control the whole of $X$, we usually need to have good information not only on general members of the family, but also on the potential degeneration of the family. Thus it may look more reasonable to require, in Problem \ref{p.VMRT},  some additional properties on the behavior of the family $\sC_x$ under degeneration.  But such additional conditions would diminish the true interest of Problem \ref{p.VMRT}. This is because in the context of intrinsic geometry of uniruled manifolds, the
     properties of $\sC_x$ we are looking for should be checkable by deformation theory of curves. Deformation theory of rational curves works well at general points of nonsingular varieties, but not so at special points.
     Thus it is important to find conditions for $\sC_x$ only for general $x \in X$ in Problem \ref{p.VMRT}.
      But then  controlling the whole of $X$ using the algebraic behavior of $\{ \sC_x \subset \BP T_x(X), \mbox{ general } x \in X \} $  becomes a serious issue.

   This was exactly the issue  puzzling me when I first encountered a version of Problem \ref{p.VMRT}  about twenty years ago. At that time,
    I was working on the deformation  rigidity of Hermitian symmetric spaces in the setting of algebraic geometry. I refer the reader to \cite{Hwang06} for the details on this rigidity  problem. Here it suffices to say that the deformation rigidity of Hermitian symmetric spaces was a question originated from Kodaira-Spencer's work in 1950's and the question itself did not involve rational curves. I was trying to attack this question employing Mori's approach of   minimal rational curves, which naturally led to a version of Problem \ref{p.VMRT} when $S$ is an irreducible Hermitian symmetric space. In the setting of this rigidity question, I could derive a certain amount of algebro-geometric information on $\sC_x$ for general $x \in X$, but I was unable to figure out how to proceed from there, essentially because of the above difficulty, that it is hard to control the whole of $X$ by algebro-geometric information on $\sC_x$ for general $x \in X$.

    There was one hope. A few years earlier, Ngaiming Mok had overcome an obstacle of  a similar kind in \cite{Mok88}. In that work, Mok solved
    what is called the generalized Frankel conjecture, which asks for a characterization of Hermitian symmetric spaces among K\"ahler manifolds in terms of a curvature condition. The Frankel conjecture itself is the K\"ahler version of the Hartshorne conjecture and was settled by Siu-Yau \cite{SiuYau} around the time Mori solved the Hartshorne conjecture. Since the method used by Siu and Yau was rather restrictive, Mok naturally took the approach of Mori and encountered a situation similar to  Problem \ref{p.VMRT}. Now in his situation, there is  a Riemannian metric on $X$ and Mok could relate $\sC_x$ to
    a suitably deformed Riemannian metric.  This enables him to show that $X$ is Hermitian symmetric space using Berger's work on Riemannian holonomy.
    Roughly speaking, in \cite{Mok88}
\begin{quote} \em    the difficulty in Problem \ref{p.VMRT}  was overcome by relating $\sC_x$ to a Riemannian structure.  \end{quote}
This shows that differential geometry can be a recourse for Problem \ref{p.VMRT} when $S$ is a Hermitian symmetric space. Indeed, compared with  tools in algebraic geometry, methods of differential geometry tend to be more effective when the available data are only at general points of a manifold.
    Motivated by this, I  tried to imitate Mok's argument in the setting of the deformation rigidity problem. However, the nature of the deformation rigidity problem is purely algebro-geometric and it is very hard to relate it to  Riemannian structures. As a matter of fact, there had already been some unsuccessful attempts in 1960's to use Riemannian structures for the deformation rigidity question.

    This was precisely the problem I was agonizing over when I attended my first ICM: Z\"urich  1994. Having come to the congress just to have fun listening to the new developments in mathematics, I found that Mok was there as a speaker and managed to
    have a chat with him. When I  told him about the above difficulty in applying the approach of \cite{Mok88}  to the deformation rigidity problem, he gave me an enlightening comment: besides the Riemannian metric,
     there is another  differential geometric structure, a certain
    holomorphic G-structure,   which can be used to characterize a Hermitian symmetric space. His suggestion was that  one might be able to construct these G-structures using the information on  $\sC_x$ for general $x \in X$ and from this to recover Hermitian symmetric spaces. This suggestion looked  promising because algebro-geometric data are closer to holomorphic structures  than
    Riemannian structures.

    Soon after the congress, I started looking into G-structures. I realized that there is a far-reaching generalization of Riemannian structures by Elie Cartan and the G-structures modeled on Hermitian symmetric spaces are special examples of Cartan's general theory of geometric structures.  To recover these G-structures, it was necessary to investigate the geometry of  $\sC_x$'s in depth in the setting of  the deformation rigidity problem.  I had subsequent
    communications with Mok, and we started working together on this problem. Our  collaboration was successful, leading to a solution of the deformation rigidity problem in \cite{HwangMok98}.  But the most exciting point in our work was not the deformation rigidity itself.  As  mentioned, the essential part of  \cite{HwangMok98} is to construct on $X$ the G-structures modeled on Hermitian symmetric spaces.  It turns out that a crucial point of this construction lies in a study of
    the behavior of $\sC_x$'s not just as a family of projective algebraic varieties, but as data imposed on the
    tangent bundle of an open subset of $X$. In other words, we had to treat these $\sC_x$'s as if
\begin{quote} \em the union of $\sC_x$'s for general $x \in X$ is a  differential geometric structure. \end{quote}
    And why not? Such a family of subvarieties in $\BP T_x(X)$ is a legitimate example of  Cartan's general geometric structures!
    So what happened can be summarized as follows. Initially we had been trying to relate $\sC_x$'s to some differential geometric structures. These differential geometric
     structures were Riemannian structures in \cite{Mok88} and then G-structures in \cite{HwangMok98}.
       But actually, they have been there all along, namely, $\sC_x$'s themselves!

      Now once we accept $\sC_x$'s as a differential geometric structure,  there is no need to restrict ourselves
     to Hermitian symmetric spaces.  This geometric structure exists for any
     uniruled projective manifold $S$ and its minimal rational curves! This was an epiphany for me.
     We realized that the variety $\sC_x$ deserves a name of its own and endowed it with the appellation, somewhat uncharming,  `variety of minimal rational tangents'.

       Realizing the varieties of minimal rational tangents as geometric structures opens up an approach to Problem \ref{p.VMRT} via Cartan geometry. In fact, Mok and I were able to show  in \cite{HwangMok01}  (see Theorem \ref{t.CartanFubini} for a precise statement)
\begin{theorem}\label{t.roughCF}
     Assume that $S$ is a fixed uniruled projective manifold with $b_2(S) =1$ and $\sC_s$ at a general point $s \in S$ is a smooth irreducible variety of positive dimension. If $X$ is a uniruled projective manifold with $b_2(X) =1$ and  the differential geometric structures defined by $\sC_s$'s and $\sC_x$'s are locally equivalent in the sense of Cartan, then  $S$ and $X$ are biregular. \end{theorem}
     This means that for a large and interesting class of uniruled projective manifolds, Problem \ref{p.VMRT} can be solved by studying  the Cartan geometry of the structures  defined by $\sC_x$'s.
      By Theorem \ref{t.roughCF}, the essence of Problem \ref{p.VMRT} has become
\begin{quote} \em searching for  algebro-geometric properties of varieties of minimal rational tangents which make it possible to control the Cartan geometry of the geometric structures defined by them.\end{quote}
     As we will see in Section 4, this search has been successful in a number of cases  and Problem \ref{p.VMRT} has been answered for  some uniruled projective manifolds, including irreducible Hermitian symmetric spaces.

     Since \cite{HwangMok98}, the theory of varieties of minimal rational tangents has seen exciting developments and  has found a wide range of applications in algebraic geometry.
     For  surveys on these developments and applications, we refer the reader to \cite{Hwang01},
     \cite{HwangMok99},  \cite{KebekusSola} and \cite{Mok08c}.  The purpose of this article
      is to give an introduction to one special aspect of the theory, the development centered around Problem \ref{p.VMRT}. This is a special aspect, because many results on varieties of minimal rational tangents and their applications are not directly related to it.
      Yet, this is the most fascinating aspect: it offers an area for  a fusion of algebraic geometry and differential geometry, more specifically, a fusion of Mori geometry of minimal rational curves and Cartan geometry of cone structures.
       We will stick to the core of this aspect and will not go into the diverse issues arising from it. Those interested in further directions of explorations may find my MSRI article \cite{Hwang12} useful.

\medskip
{\bf Conventions}
We will work over the complex numbers and all our objects are holomorphic. Open sets refer to the Euclidean topology, unless otherwise stated. All manifolds are connected. A projective manifold is a smooth irreducible projective variety. A variety is a complex analytic set which is not necessarily irreducible, but has finitely many irreducible components.
 A general point of a manifold or an irreducible variety means a point in a dense open subset.

\section{Cartan geometry: Cone structures}
{\em A priori}, this section is about local differential geometry and has nothing to do with rational curves.
We will introduce a class of geometric structures,  cone structures, and some related notions.
In a simpler form,  cone structures  have already appeared in twistor theory (see \cite{Manin}), but as they are not
widely known, I will try to give a detailed introduction.

\begin{definition}\label{d.smoothcone} For a complex manifold $M$, let $\pi: \BP T(M) \to M$ be the projectivized tangent bundle.
A {\em smooth cone structure} on $M$ is a  closed nonsingular subvariety $\sC \subset \BP T(M)$  such that all components of $\sC$ have the same dimension and
the restriction  $\varpi := \pi|_{\sC}$ is a submersion. \end{definition}

We may restrict our discussion  to smooth cone structures. Understanding the geometry of smooth cone structures is already challenging and lots of examples of smooth cone structures remain uninvestigated. To have a satisfactory general theory, however, we need to allow certain singularity in $\sC$. This necessitates  the following somewhat technical definition.
(Readers not familiar with singularities may skip this definition and just stick to Definition \ref{d.smoothcone}, regarding $\nu: \widetilde{\sC} \to \sC$
as an identity map in the subsequent discussion.)

\begin{definition}\label{d.cone} A {\em cone structure} on a complex manifold $M$ is a closed subvariety $\sC \subset \BP T(M)$ the normalization $\nu: \widetilde{\sC} \to \sC$ of which satisfies the following conditions. \begin{itemize} \item[(1)]
All components of $\widetilde{\sC}$ are smooth and have the same dimension. \item[(2)] The composition $\varpi:= \pi \circ \nu
: \widetilde{\sC} \to M$ is a submersion. In particular, the relative tangent bundle $T^{\varpi} \subset T(\widetilde{\sC})$ is a vector subbundle.  \item[(3)] There is a vector subbundle $\sT \subset T(\widetilde{\sC})$ with $T^{\varpi} \subset \sT$ and
${\rm rank} (\sT) = {\rm rank}(T^{\varpi}) +1,$  such that for any $\alpha \in \widetilde{\sC}$ and any  $v \in \sT_{\alpha} \setminus T^{\varpi}_{\alpha}$, the nonzero vector
${\rm d} \varpi_{\alpha}(v) \in T_{\varpi(\alpha)} (M)$ satisfies $$[ {\rm d} \varpi_{\alpha}(v)] = \nu (\alpha) \mbox{ as elements of } \BP T_{\varpi(\alpha)} (M).$$ \end{itemize} \end{definition}
The conditions (1) and (2) say that $\sC$ is allowed to be singular but  it becomes smooth after normalization and the natural projection  to
$M$ becomes a submersion.
 Note that on $\BP T(M)$, we have the tautological line bundle $\xi \subset \pi^* T(M)$.  The condition (3) says that the quotient line bundle $\sT / T^{\varpi}$ is naturally isomorphic to $\nu^* \xi$. Another useful interpretation of the condition (3) is
in terms of the following
\begin{definition}\label{d.Sm}
Given a cone structure $\sC \subset \BP T(M)$, let ${\rm Sm}(\sC) \subset \sC$ be the maximal dense open subset such that
$$\pi|_{{\rm Sm}(\sC)} : {\rm Sm}(\sC) \to \pi({\rm Sm}(\sC))$$ is a submersion. \end{definition}
Denote by $\sT^{\BP T(M)} \subset T(\BP T(M))$ the inverse image of the tautological bundle $\xi\subset \pi^* T(M)$ under ${\rm d} \pi: T(\BP T(M)) \to \pi^* T(M)$.
Then the condition (3) means that the vector bundle $\sT^{\BP T(M)} \cap T({\rm Sm}(\sC))$ on ${\rm Sm}(\sC)$, after pulling back to $\widetilde{\sC}$ by $\nu$,  extends to a vector subbundle of $T(\widetilde{\sC})$. From this interpretation of (3),  it is easy to see that
\begin{proposition}\label{p.smooth} A cone structure $\sC \subset \BP T(M)$ is a smooth cone structure if and only if $\sC$ is normal, i.e.,  the normalization $\nu: \widetilde{\sC} \to \sC$ is biholomorphic. \end{proposition}
All three conditions (1)-(3) for cone structures are of local nature on $M$. This implies
\begin{proposition}\label{p.localize} Given a cone structure $\sC \subset \BP T(M)$ and a connected open subset $U \subset M$, the restriction $$\sC|_U := \sC \cap \BP T(U) \ \subset  \ \BP T(U)$$ is a cone structure on the complex manifold $U$. \end{proposition}
By Proposition \ref{p.localize}, we can view a cone structure as a  geometric structure on $M$.
We are interested in Cartan geometry of cone structures. In particular,  isomorphisms in  cone structures are given by the following
\begin{definition}\label{d.equivalent} A cone structures $\sC \subset \BP T(M)$ on a complex manifold $M$ is {\em equivalent} to
a cone structure $\sC' \subset \BP T(M')$ on a complex manifold $M'$  if there exists a biholomorphic map $\varphi: M \to M'$ such that the projective bundle isomorphism $\BP {\rm d} \varphi: \BP T(M) \to \BP T(M')$
 induced by the differential ${\rm d} \varphi: T(M) \to T(M')$ of $\varphi$ satisfies $\BP {\rm d} \varphi (\sC) = \sC'$. \end{definition}
 It is convenient to have a localized version of this:
 \begin{definition}\label{d.localequivalence} For a cone structure $\sC \subset \BP T(M)$ (resp. $\sC' \subset \BP T(M')$) and a
 point $x \in M$ (resp. $x' \in M'$), we say that $\sC$ {\em at} $x$ is {\em equivalent} to $\sC'$ {\em at} $x'$ if there exists a neighborhood $U \subset M$ of $x$ and a neighborhood $U' \subset M'$ of $x'$ such that the restriction $\sC|_U$ is equivalent to  $\sC'|_{U'}$ as cone structures.  We say that $\sC$ is {\em locally equivalent} to $\sC'$ if there are points $ x \in M$ and $x' \in M'$ such that $\sC$ at
 $x$ is  equivalent to $\sC'$ at  $x' $.
  \end{definition}

Let us give one simple example of a cone structure.
  Let $V$ be a vector space and $Z \subset \BP V$ be a projective variety all components of which have the same dimension such that the normalization $\widetilde{Z}$ is nonsingular. Via the canonical isomorphism $T(V) = V \times V$, the projectivized tangent bundle $\BP T(V) = V \times \BP V$ contains the subvariety $\sC := V \times Z \subset \BP T(V).$ This is a cone structure.
   Indeed the normalization $\widetilde{\sC}$ is just $V \times \widetilde{Z}$ which is smooth  and $\varpi: \widetilde{\sC} \to V$ is just the projection
   $V \times \widetilde{Z} \to V$ which is a submersion, verifying the conditions (1) and (2) of Definition \ref{d.cone}.  The tautological line bundle of $Z \subset \BP V$ induces a line bundle $\chi$ in $T(V \times \widetilde{Z})$ via the normalization
   morphism $\widetilde{Z} \to Z$ and  the subbundle $\sT = T^{\varpi} + \chi$ of $T(\widetilde{\sC})$ satisfies the condition (3).

   \begin{definition}\label{d.flat} The cone structure $V \times Z \subset \BP T(V)$ on $V$ defined above is called the {\em flat cone structure with a fiber } $Z \subset \BP V$. We will denote it by ${\rm Flat}_V^Z \subset \BP T(V)$. A cone structure on a complex manifold $M$ is {\em locally flat} if it is locally equivalent to ${\rm Flat}^Z_V $ for some $Z \subset \BP V$ with $\dim V = \dim M$. \end{definition}
   \begin{definition}\label{d.isotrivial} Let $Z \subset \BP V$ be a projective variety. A cone structure $\sC \subset \BP T(M)$ is $Z$-{\em isotrivial} if for a general $x \in M$, the fiber
   $$\sC_x = \sC \cap \BP T_x(M) \ \subset \BP T_x(M)$$ is isomorphic  to $Z \subset \BP V$ as a projective variety, i.e.,
   a suitable linear isomorphism $T_x(M) \to V$ sends $\sC_x$ to $Z$. A cone structure is {\em isotrivial} if it is $Z$-isotrivial for some $Z$. \end{definition}
   A locally flat cone structure is isotrivial. But an isotrivial cone structure needs not be locally flat.
     Some isotrivial smooth cone structures are very familiar objects in differential geometry.    When $Z \subset \BP V$ is a linear subspace of dimension $p$, a $Z$-isotrivial cone structure on $M$ is just a Pfaffian system of rank $p+1$ on $M$. It is locally flat if and only if the Pfaffian system is involutive, i.e., it comes from a foliation. When $Z \subset \BP V$ is a nonsingular
   quadric hypersurface, a $Z$-isotrivial cone structure is a conformal structure on $M$. It is locally flat if and only if
   it is locally conformally  flat. A natural generalization of the conformal structure is the cone structure modeled on an irreducible Hermitian symmetric space $S=G/P$. The isotropy action of $P$ on the tangent space $T_o(S)$ at the base point $o \in S$ has
   a unique closed orbit $\sC_o \subset \BP T_o(S).$ A $Z$-isotrivial cone structure where $Z \subset \BP V$ is isomorphic to
   $\sC_o \subset \BP T_o(S)$ is called an {\em almost}  $S$-{\em structure}.  A conformal structure is exactly an almost $S$-structure where $S$ is a nonsingular quadric hypersurface, equivalently, an irreducible Hermitian symmetric space of type IV.   The natural almost $S$-structure  $\sC \subset \BP T(S)$ given by the translate of $\sC_o$ by $G$-action  is locally flat, which can be seen by Harish-Chandra
   coordinates of irreducible Hermitian symmetric spaces (see Section (1.2) in \cite{Mok08b} for a presentation in terms of explicit coordinates for Grassmannians).
The G-structure on an irreducible Hermitian symmetric space $S$ referred to  in Section 1 is essentially equal to the cone structure
$\sC \subset \BP T(S)$.

   \medskip
   How do we check the local equivalence of two cone structures?  A general method of checking equivalence of geometric structures
   has been formulated by Elie Cartan \cite{Cartan}.  The fundamental apparatus in Cartan's method is a coframe.

   \begin{definition}\label{d.coframe} Let $V$ be a vector space and let $M$ be a complex manifold with $\dim V = \dim M$. A {\em coframe} on $M$ is a trivialization $\omega: T(M) \to M \times V$, equivalently, a $V$-valued 1-form on $M$ such that $\omega_x: T_x(M) \to V$ is an isomorphism for each $x \in M$.  We will denote by $\BP \omega: \BP T(M) \to M \times \BP V$ the trivialization of the projectivized tangent bundle induced by $\omega$. Given a coframe, there exists a ${\rm Hom}(\wedge^2V, V)$-valued function $\sigma^{\omega}$ on $M$, called the {\em structure function} of $\omega$, such that $${\rm d} \omega = \sigma^{\omega}_{\bullet}(\omega \wedge \omega).$$
  A coframe is {\em closed} if ${\rm d} \omega =0$, i.e., the structure function $\sigma^{\omega}$ is identically zero.
A coframe is {\em conformally closed} if there exists a holomorphic function $f$ on an open subset $U \subset M$ such that $f \omega$ is closed on $U$.
\end{definition}

The following is a simple consequence of the Poincar\'e lemma (see Theorem 3.4 in \cite{Hwang10}).
\begin{proposition}\label{p.closed}
Let $V^{\vee} \subset {\rm Hom}(\wedge^2 V, V)$ be the natural inclusion of the dual space of $V$  given by contracting with one factor.
When $\dim M \geq 3$, a coframe  $\omega$ is conformally closed if and only if $\sigma^{\omega}$ takes values in $V^{\vee}$.
\end{proposition}

 Although Cartan's method  is applicable to the equivalence problem for arbitrary cone structures, its actual implementation
   can be challenging, depending on the type of the cone structure. For isotrivial cone structures, however, this becomes simple:

\begin{definition}\label{d.adapted} Let $\sC \subset \BP T(M)$ be a $Z$-isotrivial cone structure for a projective variety $Z \subset \BP V$. A coframe $\omega: T(M) \to M \times V$ is {\em adapted } to the cone structure if $\BP \omega (\sC) = Z$.  \end{definition}

 \begin{proposition}\label{p.adapted} An isotrivial cone structure is locally flat if and only if after restricting to an open subset, it admits a conformally closed adapted coframe. \end{proposition}

Since an isotrivial cone structure always admits an adapted coframe, we can use Proposition \ref{p.closed}
  and Proposition \ref{p.adapted} to check the local flatness of an isotrivial cone structure.
  One difficulty here is that there may be several different adapted coframes, so we need to choose the right one.
  Different choices of adapted coframes are related by the linear automorphism group of the fiber. Let us elaborate this point.

 For a projective variety $Z \subset \BP V$, let $\widehat{Z} \subset V$ be its homogeneous cone. Denote by ${\rm Aut}(\widehat{Z}) \subset {\rm GL}(V)$ the linear automorphism group of $\widehat{Z}$ and by ${\rm aut}(\widehat{Z}) \subset {\rm gl}(V)$ its Lie algebra. Since $\widehat{Z} \subset V$ is a cone, the Lie algebra ${\rm aut}(\widehat{Z})$ always contains the
 scalars $\C$.
 When ${\rm aut}(\widehat{Z}) =\C$, a $Z$-isotrivial cone structure has a unique adapted coframe   up to multiplication by functions.  Consequently, the method of Proposition \ref{p.adapted}  essentially determines the local flatness of $Z$-isotrivial cone structures when ${\rm aut}(\widehat{Z}) =\C$.

 When ${\rm aut}(\widehat{Z}) \neq \C$,  however,  compositions with ${\rm Aut}(\widehat{Z})$-valued
functions give rise to many different choices of adapted coframes for a $Z$-isotrivial cone structure. In this case, Proposition \ref{p.adapted} is not decisive and we have to consider the problem of choosing the right coframe.  This leads to the equivalence problem for G-structures where G corresponds to the group ${\rm Aut}(\widehat{Z}) \subset {\rm GL}(V)$.  The general theory of G-structures has been developed by many mathematicians.   In particular, for the G-structures modeled on Hermitian symmetric spaces, \cite{Guillemin} and \cite{Ochiai} provide a calculable criterion for local flatness in terms of the vanishing of
certain `curvature tensors', which are more elaborate version of the structure functions $\sigma^{\omega}$.

\medskip
It turns out that the cone structures we are interested in are equipped with some additional structures.

\begin{definition}\label{d.connection}
Let $\sC \subset \BP T(M)$ be a cone structure. From the condition (3) in Definition \ref{d.cone},  we have an exact sequence of vector bundles on $\widetilde{\sC}$
$$ 0 \rightarrow T^{\varpi} \rightarrow \sT \rightarrow \nu^* \xi \rightarrow 0.$$ A line subbundle $\sF \subset \sT$ is called a {\em connection} of the cone structure if $\sF$ splits this exact sequence.  Thus  a connection exists if and only if this exact sequence splits.  \end{definition}

All the cone structures we are to meet have certain canonically defined connections. These connections will have some special properties.

\begin{definition}\label{d.characteristic}
In Definition \ref{d.connection},  $\sF$ is a {\em characteristic connection} if $[\sF, [\sF, \sT]] \subset [\sF, \sT]$ at
general points of $\widetilde{\sC}$. The inclusion means that  for any local section $f$ of $\sF$ and any local section $v$ of $T^{\varpi}$ regarded as local vector fields in some open subset of $\widetilde{\sC}$, the Lie bracket $[f, [f, v]]$ is a local section of $[\sF, \sT]$.
\end{definition}

 The most important property of a characteristic connection is its uniqueness for a large class of cone structures. This condition is formulated in terms of the Gauss map and the projective second fundamental form. Let us recall the definition.
\begin{definition} Let $Z \subset \BP V$ be an irreducible projective variety of dimension $p$. The {\em Gauss map} of $Z$ is the morphism $\gamma: {\rm Sm}(Z) \to {\rm Gr}(p+1, V)$  defined on the smooth locus of $Z$ by associating to a smooth point $\alpha$ of $Z$  the affine tangent space $T_{\alpha}(\widehat{Z}) \subset V$,  the tangent space of the homogeneous cone $\widehat{Z} \subset V$ along $\alpha.$ We say that the Gauss map of $Z$ is {\em nondegenerate} if $\gamma$ is generically finite over its image.
   Let $\alpha \in {\rm Sm}(Z)$ be a smooth point of $Z$
and let $N_{Z, \alpha}$ be the normal space of $Z$ inside $\BP V$ at $\alpha$. The differential of $\gamma$ defines a homomorphism
$${\rm II}_{Z,\alpha} : {\rm Sym}^2 T_{\alpha}(Z) \to N_{Z, \alpha},$$ called the {\em projective second fundamental form}.
We say that ${\rm II}_{Z, \alpha}$ is {\em  nondegenerate} if its null space
$${\rm Null}_{{\rm II}_{Z, \alpha}} = \{ v \in T_{\alpha}(Z), \ {\rm II}_{N_{Z, \alpha}} (v, u) = 0 \mbox{ for all } u \in T_{\alpha}(Z)\}$$ is zero. Then the Gauss map of $Z$ is nondegenerate if the projective second fundamental form of $Z$ is nondegenerate at a general point of $Z$.  It is well-known (e.g. Theorem 3.4.2 in \cite{IveyLandsberg}) that  if an irreducible projective variety $Z$ is smooth and not a linear subspace of $\BP V$, then the Gauss map of $Z$ is nondegenerate, or equivalently, the projective second fundamental form ${\rm II}_{Z, \alpha}$ is nondegenerate at a general point $\alpha \in Z$.
\end{definition}

The uniqueness result for a characteristic connection is the following result from \cite{HwangMok04}.
\begin{theorem}\label{t.unique}
Let $\sC \subset \BP T(M)$ be a cone structure such that all components of the fiber $\sC_x$ for a general $x \in M$  have nondegenerate Gauss maps. Then $\sC$ has at most one characteristic connection. \end{theorem}

It is easy to see that the flat cone structure ${\rm Flat}_V^Z$ in Definition \ref{d.flat} has a characteristic connection given by the intersection of $\sT$ with the fibers of the projection map $V \times \widetilde{Z} \to \widetilde{Z}$.
Cone structures admitting characteristic connections have certain amount of flatness, although this is not easy to explicate.  One manifestation is  the following proposition (see Theorem 6.2 in \cite{Hwang12} for a proof).  Although it is stated here without any regard to
minimal rational curves, a version of this proposition  is first discovered in \cite{HwangMok98} for varieties of minimal rational tangents and  has been  the key revelation on the significance of  the differential geometric interpretation of the varieties of minimal rational tangents, as mentioned in Section 1.

\begin{theorem}\label{t.Pfaff} Given a cone structure $\sC \subset \BP T(M)$ admitting a characteristic connection, denote by ${\rm Pf}(\sC)$ the Pfaffian system
defined on a dense open subset of $M$ by the linear span of the homogeneous cone $\widehat{\sC} \subset T(M)$. Then for any
$\alpha \in {\rm Sm}(\widehat{\sC})$ and $\beta \in T_{\alpha}(\widehat{\sC}) \cap T_{\alpha}^{\pi}$, and any local sections
 $\overrightarrow{\alpha}$ and $\overrightarrow{\beta}$ of ${\rm Pf}(\sC)$ extending $\alpha$ and $\beta$, the Lie bracket $[ \overrightarrow{\alpha}, \overrightarrow{\beta}] $ belongs to ${\rm Pf}(\sC)$ at the point $\pi(\alpha)$. \end{theorem}

 When $\sC \subset \BP T(M)$ itself is a Pfaffian system, i.e., when ${\rm Pf}(\sC) = \sC$, Proposition \ref{t.Pfaff} says that the existence of a characteristic connection on $\sC$ implies that $\sC$ is involutive. This is an example of the statement that characteristic connections contain certain amount of flatness. Another example of this phenomena is the next result from
 \cite{Hwang10}:

 \begin{theorem}\label{t.hypersurface}
 Let $Z \subset \BP V$ be a smooth hypersurface of degree $\geq 4$. Let $\sC \subset \BP T(M)$ be a $Z$-isotrivial cone structure.
 If $\sC$ has a characteristic connection, then it is locally flat. \end{theorem}

The above results show that the existence of characteristic connections impose  severe restrictions on isotrivial cone structures.  We expect similar restrictions on non-isotrivial cone structures, although no specific results are known.

\medskip
Being a characteristic connection is a local property of a connection $\sF$ on $\widetilde{\sC}$: the condition in Definition \ref{d.characteristic} is to be checked
on an open subset of $\widetilde{\sC}$. The connections we are interested in have  another  important property which is of  a global nature. To introduce this property, we note that there is a natural vector subbundle  $\sP \subset T({\rm Sm}(\sC))$ defined as follows.
At a smooth point $\alpha \in \sC$ with $x = \pi(\alpha)$, we have the differential ${\rm d} \varpi_{\alpha}: T_{\alpha}(\sC) \to T_x(M)$ of the projection $\varpi= \pi|_{{\rm Sm}(\sC)}:  {\rm Sm}(\sC) \to M$. Since the fiber $\sC_x \subset \BP T_x(M)$ of $\pi|_{\sC}$ is smooth at $\alpha$
 by Definition \ref{d.cone} (2), we have the affine tangent space $T_{\alpha}(\widehat{\sC_x}) \subset T_x(M)$. Define $\sP_{\alpha} \subset T_{\alpha}(\sC)$ by
$$\sP_{\alpha} := {\rm d} \varpi_{\alpha}^{-1}(T_{\alpha}(\widehat{\sC}_x)).$$  This defines a vector bundle $\sP$ on ${\rm Sm}(\sC).$

\begin{definition}\label{d.sP}
View  $T^{\varpi} \subset \sT$ on $\widetilde{\sC}$ and a connection $\sF \subset \sT$ as vector bundles on ${\rm Sm}(\sC)$ via the normalization morphism $\nu: \widetilde{\sC} \to \sC$ which is an isomorphism over ${\rm Sm}(\sC)$. Then we have
$\sT \subset \sP \subset T({\rm Sm}(\sC))$.
  A connection $\sF \subset \sT$ is $\sP$-{\em splitting} if there exists a vector subbundle $\sW \subset \sP$ on ${\rm Sm}(\sC)$ that splits $$0 \to T^{\varpi} \to \sP \to \sP/T^{\varpi} \to 0$$ such that $\sW \cap \sT = \sF$ on ${\rm Sm}(\sC)$.
\end{definition}
The connections we are interested in are $\sP$-splitting. The significance of the $\sP$-splitting property has been noticed only recently
in \cite{Hwang13}, and  many of its implications are yet to be discovered. It is used in \cite{Hwang13} in the following way.
\begin{theorem}\label{t.hypersurface2} Let $\sC \subset \BP T(M)$ be a smooth cone structure of codimension 1. In other words, $\sC$ is a smooth hypersurface in $\BP T(M)$ such that $\varpi: \sC \to M$ is a submersion. If $\sC$ has a $\sP$-splitting connection and $\dim M \geq 4$, then it is isotrivial. \end{theorem}
In particular, if the degree of the fiber $\sC_x$ in Theorem \ref{t.hypersurface2} is at least 4 and the $\sP$-splitting connection is also a characteristic connection, then $\sC$ is locally flat by Theorem \ref{t.hypersurface}. Actually, the  requirement in
Theorem \ref{t.hypersurface} that the degree $d$ of the hypersurface is at least 4 can be weakened to $d \geq 3$ if the connection is
$\sP$-splitting. Thus at least for smooth cone structures of codimension 1, the existence of a $\sP$-splitting connection has a significant consequence. This is to be contrasted with cone structures of codimension 1 that are not smooth. There are examples
discovered in \cite{CasagrandeDruel} of cone structures of codimension 1 that have $\sP$-splitting characteristic connections
but  are not isotrivial.

\medskip
Cartan geometry of cone structures with $\sP$-splitting characteristic connections  is our central interest from the differential geometric side.
As we will see in Section 4, there are lots of examples of cone structures with $\sP$-splitting characteristic connections. Properties of such structures are intricately related to the projective geometry of the fibers $\sC_x$. Thus this Cartan geometry has an inseparable link with projective algebraic geometry. This is analogous to the fact that Cartan geometry of G-structures has an intimate link with representation theory of Lie groups. For example, the proofs of Theorem \ref{t.hypersurface} and Theorem \ref{t.hypersurface2} use cohomological properties of smooth hypersurfaces in projective space. The number of results in
this direction is still very small and the investigation of cone structures with $\sP$-splitting characteristic connections is a wide open area.

Let us close this section with one remark. There is an additional property, which could be called the {\em admissibility} of a connection, that holds for all connections we are interested in. This property arises from Bernstein-Gindikin's admissibility condition in integral geometry \cite{BernsteinGindikin}.
 Significant consequences of admissibility have not yet been found in connection with the topic of this article, which  is the reason I have skipped discussing this property.
However, this additional condition may lead to interesting discoveries in the future.

\section{Mori geometry: minimal rational curves}\label{s.Mori}
Our major interest in algebraic geometry is in uniruled projective manifolds, i.e., projective manifolds covered by rational curves. Recall that a rational curve $C$ on a projective manifold $X$ is a curve $C \subset X$ with normalization  $\nu_C: \BP^1 \to C$ by $\BP^1$. The set ${\rm RatCurves}(X)$ of all rational curves on $X$ can be given a scheme structure and its normalization is denoted by ${\rm RatCurves}^{\rm n}(X)$. Each irreducible component $\sK$ of ${\rm RatCurves}^{\rm n}(X)$ is a quasi-projective variety equipped with the universal $\BP^1$-bundle $\rho_{\sK}: {\rm Univ}_{\sK} \to \sK$ and the associated cycle morphism
$\mu_{\sK}: {\rm Univ}_{\sK} \to X$. This means that for each $z \in \sK$, the corresponding rational curve $C \subset X$ is given by
$\mu_{\sK}(\rho_{\sK}^{-1}(z))$ and the morphism $$\nu_C:= \mu_{\sK}|_{\rho_{\sK}^{-1}(z)} : \BP^1 \to C= \mu_{\sK}(\rho_{\sK}^{-1}(z))$$ is the normalization of $C$. For a rigorous presentation of this foundational material, we refer the reader to \cite{Kollar}.

Now I am going to introduce a number of terms related to uniruled manifolds and rational curves. I should warn the reader that most of
these are {\em not standard}: they appear under different names in the literature. As is the case in any growing area of mathematics,
the technical terms have not yet been completely standardized. I believe the terms introduced below are shorter and more intuitive than some of the ones in use (including some in my own papers)
  for nonexperts to remember their meaning.  You may regard the definitions below as nicknames we will use in this article.
To start with, we can give a precise definition of a uniruled projective manifold in the following form.
\begin{definition}\label{d.uniruling}
An irreducible component $\sK$ of ${\rm RatCurves}^{\rm n}(X)$ is called a {\em uniruling} on $X$ if the cycle morphism $\mu_{\sK}: {\rm Univ}_{\sK} \to X$ is dominant. A projective manifold $X$ is {\em uniruled} if it has a uniruling. For a line bundle $L$ on $X$, we will denote by $\deg_L(\sK)$  the $L$-degree of a member of $\sK$. \end{definition}

A fundamental tool in the study of unirulings on $X$ is the deformation theory of rational curves on $X$, or equivalently, the deformation theory of morphisms $\BP^1 \to X$. By the classical Kodaira theory, given  a rational curve $C \subset X$, the first-order deformation of the normalization morphism
$\nu_C: \BP^1 \to X$ regarded as a map to $X$ is controlled by the pull-back $\nu_C^*T(X)$ of the tangent bundle of $X$.
In this regard, the following definition is fundamental.
\begin{definition}\label{d.free}
A rational curve $C \subset X$ is {\em free} if $\nu_C^*T(X)$ is semi-positive, i.e., of the form $\sO(a_1) \oplus \cdots \oplus \sO(a_n), n = \dim X,$ with $a_i \geq 0$ for all $i$. \end{definition}
Free rational curves have a nice deformation theory because  $H^1(\BP^1, \nu_C^*T(X))=0$ by the semi-positivity of $\nu_C^*T(X)$. This cohomology group contains  the obstruction to realizing deformations of $\nu_C$ from its infinitesimal deformations in $H^0(\BP^1, \nu_C^* T(X))$. Thus the vanishing implies the following.
\begin{theorem}\label{t.free} Let $\sK$ be an irreducible component of ${\rm RatCurves}^{\rm n}(X)$. Denote by $\sK^{\rm free} \subset \sK$ the parameter space of members of $\sK$ that are free.  Then $\sK$ is a uniruling if and only if
 $\sK^{\rm free}$ is nonempty. In this case, $\sK^{\rm free}$ is a Zariski open subset of the smooth locus of $\sK$. \end{theorem}
 Given a uniruling $\sK$ on $X$ and a point $x \in X$, let $\sK_x$ be the normalization of the subvariety of $\sK$ parametrizing members of $\sK$ passing through $x$.  When $x$ is a general point of $X$, the structure of $\sK_x$ is particularly nice:
 \begin{theorem}\label{t.K_x}
 For a uniruling $\sK$ on a projective manifold $X$ and a general point $x \in X$, all members of $\sK_x$ belongs to $\sK^{\rm free}$. Furthermore,  the variety $\sK_x$ is  a finite union of smooth quasi-projective varieties of dimension $\deg_{K^{-1}_X}(\sK) -2.$ \end{theorem}
 Both Theorem \ref{t.free} and Theorem \ref{t.K_x} must have been known before \cite{Mori}, although their significance has not been fully recognized until Mori's work. Now we are ready to introduce minimal rational curves.
\begin{definition}\label{d.unbreakable}
A uniruling $\sK$ on a projective manifold $X$ is {\em unbreakable } if $\sK_x$ is projective for a general $x\in X$. In other words,
$\sK$  is an unbreakable uniruling if a general fiber of the cycle morphism $\mu_{\sK}: {\rm Univ}_{\sK} \to X$ is nonempty and complete. Members of an unbreakable uniruling on $X$ will be called {\em minimal rational curves} on $X$. \end{definition}
Unbreakable unirulings exist on any uniruled projective manifold. To see this, we need the following notion.
\begin{definition}\label{d.minimal}
Let $L$ be an ample line bundle on a projective manifold $X$. A uniruling $\sK$ is a {\em minimal  with respect to} $L$,
if  $\deg_{L}(\sK)$ is minimal among all unirulings of $X$.  A uniruling is a {\em minimal uniruling} if it is minimal with respect to
some ample line bundle. Minimal unirulings exist on any uniruled projective manifold and they are unbreakable. \end{definition}
   It is essential to understand the geometric idea behind the unbreakability  of minimal unirulings.   Suppose for a  uniruling $\sK$, which is minimal with respect to an ample line bundle $L$, the variety $\sK_x$ is not projective for a general point $x \in X$. Then the members of $\sK_x$ degenerate to reducible curves all components of which are   rational curves of smaller $L$-degree than the members of $\sK$ and some components of which pass through $x$. Collecting those components passing through $x$ as $x $ varies over the
general points of $X$ gives rise to another uniruling $\sK'$ satisfying $\deg_L(\sK') < \deg_L(\sK)$, a contradiction to the minimality of $\deg_L(\sK)$. This argument gives an intuitive picture behind the definition of an unbreakable uniruling: if a uniruling is not unbreakable, its members can
be broken into members of another uniruling. More figuratively speaking, {\em if a uniruling is not unbreakable, it can be broken into
a smaller uniruling.}

It is worth introducing a special class of minimal unirulings, which are particularly interesting from the viewpoint of (extrinsic) projective algebraic geometry:
\begin{definition}\label{d.lines}
Let $L$ be an ample line bundle  on a projective manifold $X$ that is base-point free. A uniruling $\sK$ on $X$ is
  a {\em uniruling by lines} if $\deg_L(\sK) =1$.
   Geometrically, this means that there is a morphism $j: X \to \BP^N$ which is finite over $j(X)$ and sends members of $\sK$  to lines in $\BP^N$.
 \end{definition}
 Most of the classical examples of unbreakable unirulings are unirulings by lines and the morphism  $j$ is often an embedding. For example,  smooth complete intersections of low degree in $\BP^N$ are covered by lines and so are  Grassmannians under the Pl\"ucker embedding. But there are  many examples of unirulings by lines where $j$ is not an embedding.
 Also there are many minimal unirulings which  are not unirulings by lines: hypersurfaces of degree $n$ in $\BP^{n}$ have
 minimal unirulings by conics. Also there are many unbreakable unirulings that are not minimal. Trivial examples can be constructed on a product $X = X_1 \times X_2$ of two uniruled manifolds. There are more interesting examples in \cite{CasagrandeDruel} and \cite{Lau} where these non-minimal  unbreakable unirulings
rather than minimal unirulings play crucial roles. We refer the reader to \cite{IskoProk} and \cite{Kollar} for many examples of unbreakable unirulings.

All these examples illustrate that unbreakable unirulings and minimal rational curves, which exist on any uniruled projective manifolds,  are genuine extensions of  the classical notion of unirulings by lines:
$$ \{ \mbox{unirulings by lines}\} \  \subset \{ \mbox{ minimal unirulings } \} \  \subset \{ \mbox{ unbrekable unirulings }\}.$$
The important point is that this extension retains many geometric properties of unirulings by lines. The most fundamental example is the following.
Let $x \neq y$ be two distinct points on a projective manifold $X$.
If $\sK$ is a uniruling by lines on $X$, we know that there exists at most one member of $\sK$ through $x$ and $y$.
Mori shows that a weaker version of this property continues to hold for unbreakable unirulings:
\begin{theorem}\label{t.BB}
Let $\sK$ be an unbreakable uniruling. Then for a general point $x \in X$ and any other point $y \in X$, there does not exist
a positive-dimensional family of members of $\sK$ that pass through both $x$ and $y$. In particular, $\dim \sK_x \leq \dim X -1$.
\end{theorem}
Theorem \ref{t.BB} is proved by what is called the `bend-and-break' argument. Geometrically, it says that any 1-dimensional
family of rational curves which share two distinct points in common  must degenerate into a reducible curve.
This  is the most important geometric property of an unbreakable uniruling.
Combined with Theorem \ref{t.K_x}, the bound on $\dim \sK_x$ in Theorem \ref{t.BB} implies that $\deg_{K^{-1}_X}(\sK)  \leq \dim X +1$ for any unbreakable uniruling. In other words, the  $\sK_X^{-1}$-degrees of minimal rational curves are bounded by $\dim X +1$. The fact that all uniruled projective manifolds are covered by rational curves of small degree is of fundamental significance and has eventually developed into the boundedness of Fano manifolds, for which we refer the reader to Chapter V of \cite{Kollar}.

 An infinitesimal version of  Theorem \ref{t.BB} is important for us.
 A key notion here is the following
\begin{definition}\label{d.unbending}
A rational curve $C \subset X$ is {\em unbending} if under the normalization $\nu_C: \BP^1 \to C \subset X$, the vector bundle
$\nu_C^* T(X)$ has the form $$\nu_C^*T(X) \cong \sO(2) \oplus \sO(1)^p \oplus \sO^{n-1-p}$$ for some integer $p$
satisfying $0 \leq p \leq n-1,$ where $n= \dim X$. \end{definition}
What is the rationale behind the name `unbending'? Note that a free rational curve $C \subset X$ is unbending if and only if
for any two distinct points $x \neq y \in \BP^1$ and their maximal ideals ${\bf m}_x$ and ${\bf m}_y$ in $\sO_{\BP^1}$,  $$H^0(\BP^1, \nu_C^*T(X) \otimes {\bf m}_x \otimes {\bf m}_y) = H^0(\BP^1, T(\BP^1) \otimes {\bf m}_x \otimes {\bf m}_y).$$
In the standard deformation theory, this means that $C$ does not have infinitesimal deformation fixing two distinct points.  Figuratively, we can say that `an unbending curve cannot be bent infinitesimally'.
Comparing this with Theorem \ref{t.BB}, we wonder whether members of an unbreakable uniruling are unbending.
This is indeed the case for general members:
\begin{theorem}\label{t.unbending}
A general member of an unbreakable uniruling is unbending. \end{theorem}
Theorem \ref{t.unbending} is not a direct consequence of Theorem \ref{t.BB} because the notion of an unbending rational curve gives
information on  infinitesimal deformation only, while Theorem \ref{t.BB} is concerned with an actually realized deformation.
Theorem \ref{t.unbending} enables us to control the behavior of the tangent directions of members of an unbreakable uniruling.  To
study this  behavior systematically, it is convenient to introduce the tangent map.
\begin{definition}\label{d.tangent}
For any uniruling $\sK$ on a projective manifold $X$ and a point $x \in X$,
the rational map $\tau_x: \sK_x \dasharrow \BP T_x(X)$ sending a member of $\sK_x$ that is smooth at $x$ to its tangent direction is called the {\em tangent map} at $x$.
\end{definition}
If $C$ is an unbending member of $\sK_x$, the tangent map can be extended to $[C] \in
\sK_x,$ even when $C$ is singular at $x$, because  the differential ${\rm d} \nu_C: T(\BP^1) \to \nu_C^*T(X)$ is injective.
In fact, a stronger result holds.
\begin{theorem}\label{t.tangent}
In Definition \ref{d.tangent}, if $C$ is an unbending member of $\sK_x$, the tangent map $\tau_x$ is well-defined and immersive at $[C] \in \sK_x$. \end{theorem}
In particular, Theorem \ref{t.unbending} implies that for an unbreakable uniruling $\sK$ and a general point $x \in X$, the tangent map $\tau_x$ is generically finite over its image.

\medskip
  It is worth comparing Theorem \ref{t.tangent} with Theorem \ref{t.free}. Theorem \ref{t.free} suggests that the freeness of rational curves is an `individualized' version of the notion of a uniruling. In a similar way, Theorem \ref{t.BB} suggests that the unbending property of rational curves is an `individualized' version of the  notion of an unbreakable uniruling. The correspondence in the
    latter case, however, is less exact:  general members of an unbreakable uniruling are unbending, but a uniruling whose general members are unbending is not necessarily unbreakable.

    There is another important difference between the two correspondences.
      When $\sK$ is a uniruling, we know that every member of $\sK_x$ for a general $x \in X$ is free by Theorem \ref{t.K_x}.
      Is it true that for an unbreakable uniruling $\sK$, every member of $\sK_x$ for a general $x \in X$ is unbending?
    Using Theorem \ref{t.tangent}, we can formulate this question as follows:
    \begin{question}\label{Q} For an unbreakable uniruling $\sK$ on a projective manifold $X$, is the tangent map $\tau_x: \sK_x \dasharrow \BP T_x(X)$ at a general point $x\in X$ extendable to an immersion $\tau_x: \sK_x \to \BP T_x(X)$? \end{question}
      It is easy to see that if $\sK$ is a uniruling by lines, then the answer is yes: $\tau_x$ is an embedding. This corresponds to the
      classical property of lines that they are determined by their tangent directions at a given point.
       Encouraged by this special case, it has been expected
      that the answer to Question \ref{Q} is affirmative for all unbreakable unirulings, or at least for all minimal unirulings. Recently, however,  counterexamples have
      been discovered: an unbreakable one in \cite{CasagrandeDruel}  and then a minimal one in \cite{HwangKim14}.

\medskip
 Although not all members of $\sK_x$ are as nice as we would wish them to be, they are still considerably well behaved, as  the following result of Kebekus shows. In \cite{Kebekus}, an in-depth analysis of singularities of members of $\sK_x$ has been carried out. Among other things, Kebekus has shown
 \begin{theorem}\label{t.Kebekus}
 For an unbreakable uniruling $\sK$ and a general point $x \in X$, all members of $\sK_x$ are immersed at the point corresponding to $x$. \end{theorem}
 To elaborate, a member of $\sK_x$ is given by a morphism $\nu_C: \BP^1 \to C \subset X$ with a point $o\in \BP^1$ satisfying $\nu_C(o) = x$. Theorem \ref{t.Kebekus} says that $({\rm d} \nu_C)_o: T_o(\BP^1) \to T_x(X)$ is injective.  Using this, Kebekus
 has shown the following important result. \begin{theorem}\label{t.Kebekus2}
 In the setting of Theorem \ref{t.Kebekus}, the tangent morphism $\tau_x: \sK_x \to \BP T_x(X)$ can be defined by assigning to each member $C$ of $\sK_x$ its tangent direction $$\BP ({\rm d} \nu_C (T_o(\BP^1))) \in \BP T_x(X).$$
This morphism $\tau_x$ is finite over its image. \end{theorem}

 Kebekus's study of singularities of members of $\sK_x$, including Theorem \ref{t.Kebekus}, plays an important role in the proof of Theorem \ref{t.CMSB} due to Cho, Miyaoka and Shepherd-Barron \cite{CMSB}. In terms of Theorem
 \ref{t.Kebekus2}, we can state it as follows.
 \begin{theorem}\label{t.CMSB2}
 Let $X$ be a uniruled projective manifold and let $\sK$ be an unbreakable uniruling. Suppose for a general point $x \in X$,
 the tangent morphism $\tau_x: \sK_x \to \BP T_x(X)$ is dominant. Then $X$ is projective space and $\sK$ is the space of lines. \end{theorem}
 By Theorem \ref{t.Kebekus2}, the morphism $\tau_x$ in Theorem \ref{t.CMSB2} is a finite covering of projective space. The essential point in the proof of Theorem \ref{t.CMSB2} is to prove that $\tau_x$ is a birational morphism, thus an isomorphism.
 In fact, it is fairly easy to show that $X$ is projective space if $\tau_x$ is an isomorphism. In this sense,
 the following result in \cite{HwangMok04} is a generalization of Theorem \ref{t.CMSB2}.
 \begin{theorem}\label{t.birational}
 Let $X$ be a uniruled projective manifold and let   $\sK$ be an unbreakable uniruling. For a general point $x \in X$, the tangent
 morphism $\tau_x: \sK_x \to \BP T_x(X)$ is birational over its image. \end{theorem}
 Combining Theorem \ref{t.Kebekus2} and Theorem \ref{t.birational}, we see that $\tau_x$ is the normalization of its image
 in $\BP T_x(X)$.

 The collection of tangent morphisms $\{ \tau_x, \mbox{ general } x \in X\}$ can be assembled into a single map. Recall that we have
 the universal $\BP^1$-bundle $\rho_{\sK}: {\rm Univ}_{\sK} \to \sK$ and the cycle morphism $\mu_{\sK}: {\rm Univ}_{\sK} \to X$. When
 $\sK$ is unbreakable, the variety
   $\sK_x$ for a general $x \in X$ can be identified with $\mu_{\sK}^{-1}(x) \subset {\rm Univ}_{\sK}$ and we have a rational map $\tau: {\rm Univ}_{\sK} \dasharrow \BP T(X)$ with a commuting diagram
   $$ \begin{array}{ccc} {\rm Univ}_{\sK} & \stackrel{\tau}{\dasharrow} & \BP T(X) \\
   \mu_{\sK} \downarrow & & \downarrow \pi \\ X & = & X \end{array} $$
   such that the tangent morphism $\tau_x: \sK_x \to \BP T_x(X)$ is just the fiber of this diagram at a general point $x \in X$. In terms of $\tau$, we can summarize Theorem \ref{t.K_x},  Theorem \ref{t.Kebekus2} and Theorem \ref{t.birational} as follows.
 \begin{theorem}\label{t.summary}
 Let $\sK$ be an unbreakable uniruling on a projective manifold $X$ and let $\tau: {\rm Univ}_{\sK} \dasharrow  \BP T(X)$ be the tangent map.
 Then there exists a Zariski open subset $X_o \subset X$ satisfying the following properties. \begin{itemize} \item[(i)]
 The dense open subset ${\rm Univ}_{\sK}^{o} := \mu_{\sK}^{-1}(X_o) $ of ${\rm Univ}_{\sK}$  is a smooth quasi-projective variety.  \item[(ii)] The morphism $\mu_{\sK}$ restricted to ${\rm Univ}_{\sK}^{o}$ is a submersion over $X_o$. \item[(iii)] The restriction of
 $\tau$ gives a morphism  $\tau^o: {\rm Univ}_{\sK}^{o} \to \BP T(X_o)$ which is normalization of its image.
 \item[(iv)] The two vector subbundles $T^{\mu_{\sK}}$ and $T^{\rho_{\sK}}$ of $T({\rm Univ}_{\sK}^o)$ are transversal.  \end{itemize}
\end{theorem}
In fact, (i) and (ii) follow from Theorem \ref{t.K_x}, (iii) follows from Theorem \ref{t.birational} and finally (iv) follows from
Theorem \ref{t.Kebekus}.

\medskip
In this section, we have collected some of the key results on minimal rational curves that are needed for the next section.
There are many other results on minimal rational curves which we have omitted, for which we refer the reader to \cite{Kollar} and
the survey papers cited  in Section 1.

\section{From Mori to Cartan: VMRT-structures}
In Section 3, we have seen that starting from a uniruled projective manifold $X$ by choosing an unbreakable uniruling $\sK$ and a general point $x \in X$, we obtain $\sK_x$, a finite union of projective manifolds.  Since $\dim \sK_x < \dim X$ by Theorem
\ref{t.BB}, the variety $\sK_x$ is likely to be simpler than $X$. This opens up the possibility of using $\sK_x$ to study the geometry of $X$. Indeed, at least for unirulings by lines, the variety $\sK_x$ has been used in this way in many classical constructions.
Now Theorem \ref{t.Kebekus2} and Theorem \ref{t.birational} show a bigger
 advantage of this $\sK_x$: it is  provided with a morphism $\tau_x: \sK_x \to \BP T_x(X)$ which is almost an embedding, a normalization of its image.
 That is, starting from the intrinsic information of an unbreakable uniruling, we obtain a natural projective subvariety
 ${\rm Im}(\tau_x) \subset \BP T_x(X)$.
   The extrinsic projective geometry of ${\rm Im}(\tau_x)$ yields intrinsic information
 on $X$ and $\sK$. Since $\tau_x$ is just the normalization of its image, it seems more advantageous to
 look at ${\rm Im}(\tau_x)$ rather than $\sK_x$. This motivates the following definition.
 \begin{definition} In the above setting, the image ${\rm Im}(\tau_x)$ is denoted by $\sC_x \subset \BP T_x(X)$ and called
 the {\em variety of minimal rational tangents} (abbr. VMRT) at $x$ associated to $\sK$. \end{definition}
 This shift of attention from $\sK_x$ to $\sC_x$ connects  Mori geometry to Cartan geometry:
\begin{definition}\label{d.vmrt}
In the setting of Theorem \ref{t.summary},  put $\widetilde{\sC} := {\rm Univ}_{\sK}^o$ and $\sC := \tau({\rm Univ}^o_{\sK})$.
Then $\nu= \tau^o : \widetilde{\sC} \to \sC$ is a normalization morphism. Putting $\sT := T^{\rho_{\sK}} \oplus T^{\mu_{\sK}}$ on
$\widetilde{\sC}$, Theorem \ref{t.summary} (iv) says that $\sT$ is a vector subbundle of $T(\widetilde{\sC})$. It satisfies the condition (3) of Definition \ref{d.cone}, so $\sC \subset \BP T(X_o)$ is a cone structure on the complex manifold $X_o$. This cone structure is called the {\em VMRT-structure} of the unbreakable uniruling $\sK$. Moreover $T^{\rho_{\sK}}$
gives a connection $\sF$ on this cone structure, called the {\em tautological connection} on $\sC$. \end{definition}
 The following is proved in \cite{Hwang13} and \cite{HwangMok04}.
 \begin{theorem}\label{t.vmrtconnection}
 The tautological connection $\sF$ on the VMRT-structure in Definition \ref{d.vmrt} is a $\sP$-splitting characteristic connection. \end{theorem}

       It follows that a choice of an unbreakable uniruling on a uniruled projective manifold gives rise to a cone structure with
       a natural $\sP$-splitting characteristic connection: the VMRT-structure with the tautological connection. This provides  diverse examples of such cone structures. How much  diversity?  We will see shortly that  this transition from algebraic geometry to differential geometry
         $$ (X, \sK) \Rightarrow (\sC \subset \BP T(X_o), \sF)$$ is injective, under some topological assumptions.
Thus one may say that such cone structures are almost as diverse as all uniruled projective manifolds and all unbreakable unirulings on them.
Referring the reader to \cite{Hwang01} and \cite{HwangMok99}   for many interesting examples of VMRT-structures,
       let us just note  three  salient features of this diversity.

         Firstly, for any irreducible smooth projective variety $Z \subset \BP V$, there exists a  VMRT-structure
       which is $Z$-isotrivial and locally flat (Example 1.7 in \cite{Hwang10}). In fact, writing $W = V \oplus \C, $ we can  regard  $Z \subset \BP V \subset \BP W$
       as a smooth subvariety of $\BP W$ contained in the hyperplane $\BP V$. Then the blow-up $X:={\rm Bl}_Z(\BP W)$ has
       an unbreakable uniruling  whose general members are proper transforms of lines in $\BP W$ intersecting $Z$ at one point. Then for
       $x \in X$ over a point on $\BP W \setminus \BP V$, the VMRT $\sC_x \subset \BP T_x(X)$ is isomorphic to $Z \subset \BP V$
       as projective varieties. It is easy to see that this VMRT-structure is locally flat.

Secondly, there are isotrivial VMRT-structures which are not locally flat. A simple example of this type is provided by a homogeneous contact manifold $X$ with second Betti number 1, different from projective space. This is a homogeneous projective manifold $X= G/P$ equipped with a $G$-invariant holomorphic contact structure $\sD\subset T(X)$, i.e., a Pfaffian system of corank 1 which is maximally non-integrable. It has a uniruling by lines,  the VMRT-structure $\sC \subset \BP T(X)$ of which is  $G$-invariant and hence it is isotrivial. Furthermore,  at each $x \in X$,
the VMRT $\sC_x$ is contained in $\BP \sD_x$ and, in fact, spans $\BP \sD_x$ (see \cite{Hwang97}). This implies
that $\sC$ is not locally flat, because if it were, then $\sD$ would be integrable.

 Thirdly, there are many examples of non-isotrivial VMRT-structures. Actually, most of the VMRT-structures on  projective manifolds with second Betti number 1 are expected to be non-isotrivial.
A most transparent example is given by the moduli variety ${\rm SU}_C(r, d)$ of stable bundles on a curve $C$ of  genus greater than 3 with rank $r$ and fixed determinant of degree $r$. When $(r,d) =1$, this is a projective manifold with second Betti number 1. There is a minimal uniruling given by `Hecke curves', certain families of stable bundles on $C$  arising from a geometric version of the Hecke correspondence introduced in \cite{NaRa}. Its VMRT's are
iterated projective bundles on $C$ constructed from the universal bundle on $C \times {\rm SU}_C(r,d)$, which shows that
the VMRT-structure is not isotrivial (see \cite{Hwang03} and \cite{HwangRamanan} for details).

These three points illustrate the diversity of VMRT-structures. But even on these specific points, our understanding of
this diversity is very limited.

       Regarding the examples of locally flat VMRT-structures,  our construction  does not work
       if $Z$ is reducible. In fact, if $\sK_x$ is reducible,  there is a global constraint on the VMRT-structure coming from the irreducibility of $\sK$.
       For example, we cannot have $\sC_x$ consisting of two components only one of which is linear.
       In fact, if this were to happen, then collecting the linear components for  general $x \in X$, we would get an irreducible component of
       $\sC$ different from $\sC$, a contradiction to the irreducibility of $\sC$. Also the construction does not work if
       $Z$ is singular with smooth normalization. Which singular varieties can be realized as the VMRT of an unbreakable uniruling
       is a completely open question.

 As for  the examples of non-isotrivial VMRT-structures, although we expect that `most' VMRT-structures are non-isotrivial, this has  been verified only in a small number of cases  because  checking the non-isotriviality can be tricky even when $X$ is a simple well-known variety.  We do not have a general method for doing this. An approach to prove that a certain type of VMRT-structures  are not isotrivial   has been proposed in \cite{Hwang10}, but this has not been implemented in concrete examples. So far the non-isotriviality is known only for a handful of examples, and only by explicit  direct description of the VMRT's.  Besides ${\rm SU}_C(r,d)$ explained above, there are only two more types of examples we know of.  In \cite{LandsbergRobles}, Landsberg and Robles show that the uniruling by lines on   a  general smooth hypersurface of degree $d, 3 \leq d \leq n,$ in $\BP^{n+1}$ has non-isotrivial VMRT-structure. Another case is in \cite{HwangKim13} which shows the non-isotriviality of the VMRT-structure for the uniruling by lines on a double cover  of $\BP^n$  branched along a general smooth hypersurface of degree $2m, 2 \leq m \leq n-1.$

\medskip
Now that we have many cone structures with $\sP$-splitting characteristic connections arising from unbreakable unirulings on projective manifolds, it is natural to study the local equivalence problem for these structures in the sense of Cartan, which we can formulate as follows.
\begin{definition}\label{d.equivalence}
Let $X^1$ (resp. $X^2$) be a uniruled projective manifold with an unbreakable uniruling $\sK^1$ (resp. $\sK^2$).  Let $\sC^1 \subset \BP T(X^1_o)$ (resp. $\sC^2 \subset \BP T(X^2_o)$) be the corresponding VMRT-structure with the tautological connection
$\sF^1$ (resp $\sF^2$). Suppose that there exists a connected Euclidean open subset $U^1 \subset X^1_o$ (resp. $U^2 \subset X^2_o$) and a biholomorphic map $\varphi: U^1 \to U^2$ such that $\BP {\rm d} \varphi: \BP T(U^1) \to \BP T(U^2)$ sends
$\sC^1|_{U^1}$ to $\sC^2|_{U^2}$ and $\sF^1|_{U^1}$ to $\sF^2|_{U^2}$. Then we say that the VMRT-structures $\sC^1$ and $\sC^2$ are {\em locally equivalent} and $\varphi:U^1 \to U^2$ is a {\em local equivalence map} of the two VMRT-structures. \end{definition}
This is of course a natural definition from the viewpoint of Cartan geometry. But does it have significant implications in algebraic geometry? The following result from \cite{HwangMok01} shows that this is indeed so, under some topological assumptions.
\begin{theorem}\label{t.extension}
In the setting of Definition \ref{d.equivalence}, assume that $b_2(X^1) =1$ and $\dim \sK^1_x >0$ for a general $x \in X_1$, or equivalently, $\deg_{K_{X^1}^{-1}}(\sK^1) \geq 3$. Then an equivalence map $\varphi$ can be extended to
a rational map $\Phi: X^1 \dasharrow X^2$, i.e., $\varphi= \Phi|_{U^1}.$ If furthermore, $b_2(X^2) =1$, then $\Phi$ is biregular.
\end{theorem}
In other words, within the class of uniruled projective manifolds of second Betti number 1 and unirulings of anti-canonical degree
at least 3, the correspondence $$(X, \sK) \Rightarrow (\sC \subset \BP T(X_o), \sF)$$ is injective. The requirement in Theorem \ref{t.extension} that the projective manifolds have second Betti number 1 is  necessary. Indeed, it is easy to construct examples of $X^1$ with $b_2(X^1) >1$ where an analogue of
Theorem \ref{t.extension} fails. However, this condition is not a serious handicap.
In fact, projective manifolds with second Betti number 1 form a large class of projective varieties where general structure theories
of higher dimensional algebraic geometry, like the minimal model program, give little direct information.
  That VMRT-theory is effective for uniruled projective manifolds with second Betti number 1  means that it could complement these general structure theories.

 The proof of Theorem \ref{t.extension} is by analytic continuation of the map $\varphi$ along members of $\sK^1$.  This analytic continuation corresponds to  the
construction of  the developing map in Cartan geometry. Here $X^2$ is regarded as the model of the geometric structure
and the map $\varphi$ is developed to $\Phi$. The topological condition  $b_2(X^1) =1$ is used to guarantee that the analytic continuation can be carried out to cover the whole of $X^1.$ The potential multi-valuedness of the analytic continuation is taken care of by the condition $\dim \sK^1_x >0$.

 From Definition \ref{d.equivalence}, the  condition required for the map $\varphi$  in Theorem \ref{t.extension}  is
  that $\BP {\rm d} \varphi$ preserves the cone structure and also the characteristic connection. These are
  differential geometric conditions. To be useful in algebraic geometry, we need a way to find algebro-geometric conditions to guarantee them.
In most interesting cases, the requirement of preserving the connection can be replaced  by algebro-geometric condition in the following form:
\begin{theorem}\label{t.CartanFubini}
Let $X^1$ (resp. $X^2$) be a uniruled projective manifold with $b_2(X^1) =1$ (resp. $b_2(X^2) =1$). Let  $\sK^1$ (resp. $\sK^2$) be an unbreakable uniruling on $X^1$ (resp. $X^2$).  Let $\sC^1 \subset \BP T(X^1_o)$ (resp. $\sC^2 \subset \BP T(X^2_o)$) be the corresponding VMRT-structure. Suppose that there exists a connected Euclidean open subset $U^1 \subset X^1_o$ (resp. $U^2 \subset X^2_o$) and a biholomorphic map $\varphi: U^1 \to U^2$ such that $\BP {\rm d} \varphi: \BP T(U^1) \to \BP T(U^2)$ sends
$\sC^1|_{U^1}$ to $\sC^2|_{U^2}$. If the VMRT $\sC^1_x$ at a general point $x \in U^1$ has positive dimension and  nondegenerate Gauss map, then there exists a biregular morphism $\Phi: X^1 \to X^2$ such that $\Phi|_{U^1} = \varphi.$ \end{theorem}
This is a simple combination of  Theorem \ref{t.unique} and Theorem \ref{t.extension}.
Note that one Cartanian condition in Theorem \ref{t.extension}, i.e., that $\varphi$ preserves the connection $\sF$, has been replaced by an algebro-geometric condition on the fiber $\sC_x$, the nondegeneracy of Gauss map.
The latter condition holds for the VMRT of a large class of uniruled projective manifolds and unbreakable unirulings, so we may concentrate on VMRT-structures satisfying Theorem \ref{t.CartanFubini}.
Thus our main problem becomes finding algebro-geometric conditions for the local equivalence of the cone structures.
In this regard, the most central problem connecting Mori geometry and
Cartan geometry is the following:
\begin{problem}\label{p.central}
Let $X^1$ (resp. $X^2$) be a uniruled projective manifold with an unbreakable uniruling $\sK^1$ (resp. $\sK^2$).  Let $\sC^1 \subset \BP T(X^1_o)$ (resp. $\sC^2 \subset \BP T(X^2_o)$) be the corresponding VMRT-structure with the tautological connection
$\sF^1$ (resp $\sF^2$). Suppose that there exists a connected Euclidean open subset $U^1 \subset X^1_o$ (resp. $U^2 \subset X^2_o$) and a biholomorphic map $\varphi: U^1 \to U^2$ such that $  \sC^1|_{U^1}$ and $ \sC^2|_{U^2}$ are isomorphic as families of projective varieties. More precisely, suppose that we have a biholomorphic map $\psi: \BP T(U^1) \to \BP T(U^2)$  satisfying $\psi(\sC^1|_{U^1}) = \sC^2|_{U^2}$ with a commuting diagram $$\begin{array}{ccc}
\BP T(U^1) & \stackrel{\psi}{\longrightarrow} & \BP T(U^2) \\ \downarrow & & \downarrow \\ U^1 & \stackrel{\varphi}{\longrightarrow} & U^2 \end{array}.$$
Does this imply that the two cone structures
 $\sC^1|_{U^1}$ and $\sC^2|_{U^2}$ are locally equivalent? In other words, can we choose $\varphi$ such that $\psi = \BP {\rm d} \varphi$?
\end{problem}

 What makes Problem \ref{p.central} exciting and challenging is that it needs insights and techniques for the fusion of  the algebraic geometry of the projective variety $\sC_x \subset \BP T_x(X)$ and the differential geometry of the cone structure $\sC$ via the connection $\sF$. At present, no plausible uniform method to handle Problem \ref{p.central} seems conceivable and, depending on the type of the family of varieties $\sC_x$'s, different tools are required. A practicable approach is to try the cases
 where the projective geometry of the family $\sC_x$'s is simple and the differential geometric machinery of the cone structure
  is available.  A most reasonable candidate is the case when the VMRT-structure is  $Z$-isotrivial for a projective variety $Z \subset \BP V$ of simple type.

    Even for $Z$-isotrivial cases,  Problem \ref{p.central}  is highly challenging. The answer is not always affirmative. There are examples of two uniruled projective manifolds with $Z$-isotrivial VMRT-structures for the same
smooth irreducible  $Z\subset \BP V,$ which are not locally equivalent.  Since we have seen that a locally flat $Z$-isotrivial VMRT-structure exists for any given $Z \subset \BP V$,  it suffices to give examples of $Z$-isotrivial VMRT-structures that are not locally flat. We have already mentioned that homogeneous contact manifolds provide such examples where $Z \subset \BP V$ is contained in a hyperplane in $\BP V$.  There are also examples where $Z $ spans $\BP V$. A simple example is provided by  symplectic Grassmannians. Given a vector space $W$ with a symplectic form $\omega$, let us denote by ${\rm Gr}_{\omega}(k, W)$ the set of
$k$-dimensional subspaces of $W$ isotropic with respect to $\omega$. This projective manifold ${\rm Gr}_{\omega}(k,W)$ has a uniruling  by lines. The associated VMRT-structure is $Z$-isotrivial for some $Z$ because ${\rm Gr}_{\omega}(k,W)$ is a homogeneous space.  It turns out that $Z \subset \BP V$ is nondegenerate. If $2k< \dim W$, the automorphism group ${\rm Aut}(\widehat{Z})$ acts on $Z$ with two orbits  and the unique closed orbit $Z' \subset Z$ is degenerate. The corresponding subvariety $\sC' \subset \sC$ of the VMRT-structure spans a Pfaffian system $\sD \subset T({\rm Gr}_{\omega}(k,W))$ which is not integrable. This non-integrability implies that the VMRT-structure is not locally flat (see \cite{HwangMok05} for details).

 This example shows that the answer to Problem \ref{p.central} for  $Z$-isotrivial cases would depend on the variety $Z$.  This raises the following purely local problem in Cartan geometry:
 \begin{problem}\label{p.local} Find algebro-geometric conditions on a smooth irreducible $Z\subset \BP V$ which guarantee that every
$Z$-isotrivial cone structure with a $\sP$-splitting characteristic  connection is locally flat. \end{problem}

As we have already mentioned in Section 2, there are two types of $Z$ where the answer is known. When $Z$ is a linear subspace of
$\BP V$, we have seen that the corresponding Pfaffian system is involutive, thus the cone structure is locally flat.
  Another case is when  $Z \subset \BP V$ is a smooth hypersurface of degree $\geq 3$ by  Theorem \ref{t.hypersurface} together with  the comment after Theorem \ref{t.hypersurface2}.   At the moment, these are the only cases where Problem \ref{p.local} has been answered,
  even though it is reasonable to expect that it is locally flat for most choices of  $Z \subset \BP V$.

 One important case of $Z$-isotrivial VMRT-structures where an affirmative answer is obtained for Problem \ref{p.central} is when $Z$ is the VMRT of an irreducible Hermitian symmetric space. This is not done  through the local
approach of Problem \ref{p.local}. Instead,  a mixture of local and global methods
has settled the question successfully. As mentioned in Section 2, differential geometric
machinery of such structures has been developed by differential geometers and there are natural curvature tensors whose vanishing implies the local flatness.   In  \cite{Mok08a}, Mok shows that if a VMRT-structure is $Z$-isotrivial where $Z$ is the VMRT of an irreducible Hermitian symmetric space, then the $Z$-isotrivial cone structure can be extended to a neighborhood of
a general member of the unbreakable uniruling. Then the unbending property of the rational curve can be used to show the vanishing of
the curvature tensors. This way, Mok proves
\begin{theorem}\label{t.Mok08}
Let $X$ be a uniruled projective manifold with second Betti number 1. Suppose there exists an unbreakable uniruling whose
VMRT-structure is $Z$-isotrivial where $Z$ is the VMRT of an irreducible Hermitian symmetric space $G/P$. Then $X$ is biregular
to $G/P$. \end{theorem}
This resolves Problem \ref{p.VMRT} for irreducible Hermitian symmetric spaces. Mok's method has been extended to
any rational homogeneous space $G/P$ where $P$ is a maximal parabolic subgroup of a complex simple Lie group $G$, associated
to a long root of $G$ in \cite{HongHwang} and \cite{Mok08a}. Whether this remains true when $P$ is associated to a short
root of $G$, such as symplectic Grassmannians, is a most tantalizing problem in this direction. Also, it would be desirable to
recover Theorem \ref{t.Mok08} by answering Problem \ref{p.local} for the corresponding $Z$.

Combining  Theorem \ref{t.Mok08} with the affirmative answer to Problem \ref{p.local} when $Z$ is a smooth hypersurface
of degree $\geq 3$, we have the following result in \cite{Hwang13}.
\begin{theorem}\label{t.Hwang13} Let $X$ be a uniruled projective manifold of dimension $\geq 4$ with second Betti number 1. Assume that there is
an unbreakable uniruling such that the VMRT at a general point is a smooth hypersurface. Then $X$ is a quadric hypersurface. \end{theorem}

These few examples are all the results we  currently have for Problem \ref{p.central} for $Z$-isotrivial cases. They correspond to the simplest type of projective varieties, hypersurfaces and some homogeneous varieties.
Any other projective variety $Z \subset \BP V$ seems to be a new challenge. Even less -essentially no concrete result-  is known for non-isotrivial cases of
Problem \ref{p.central}.

\medskip
There is a natural extension of the equivalence problem for VMRT-structures to the setting
 of submanifold geometry.  In the investigation of any type of  geometric structures, it is natural to study submanifolds inheriting such structures. From this perspective, among submanifolds in  a uniruled projective manifold $X$, of particular importance are  uniruled submanifolds with  VMRT-structures   which are compatible with  the VMRT-structure of the ambient manifold.
The study of such submanifolds from an intrinsic viewpoint, which has been initiated in \cite{Mok08b}, is based on the fact that the analytic continuation Theorem \ref{t.extension} can be easily extended to submanifold geometry in the following sense.
\begin{theorem}\label{t.submanifold}
Let $X^1$ (resp. $X^2$) be a uniruled projective manifold with an unbreakable uniruling $\sK^1$ (resp. $\sK^2$).  Let $\sC^1 \subset \BP T(X^1_o)$ (resp. $\sC^2 \subset \BP T(X^2_o)$) be the corresponding VMRT-structure with the tautological connection
$\sF^1$ (resp $\sF^2$). Suppose that there exists a connected Euclidean open subset $U^1 \subset X^1_o$ (resp. $U^2 \subset X^2_o$) and an embedding $\varphi: U^1 \to U^2$ such that $\BP {\rm d} \varphi: \BP T(U^1) \to \BP T(U^2)$ sends
$\sC^1|_{U^1}$ into $\sC^2|_{U^2}$ and $\sF^1$ into $\sF^2$.
Assume that $b_2(X^1) =1$ and $\dim \sK^1_x >0$ for a general $x \in X_1$. Then $\varphi$ can be extended to
a rational map $\Phi: X^1 \dasharrow X^2$, i.e., $\varphi= \Phi|_{U^1}.$
\end{theorem}
The first task is to replace the preservation of connections by algebro-geometric conditions, just as Theorem \ref{t.CartanFubini}
is an improvement over  Theorem \ref{t.extension}.
 This is done in \cite{HongMok}  and the following analogue of Theorem \ref{t.CartanFubini} in submanifold geometry is obtained.
\begin{theorem}\label{t.CFi}
Let $X^1$ (resp. $X^2$) be a uniruled projective manifold  with a VMRT-structure $\sC^1 \subset \BP T(X_o^1)$ (resp. $\sC^2 \subset \BP T(X_o^2)$).   Assume that
 $b_2(X^1) =1$  and   $\dim \sC^1_x >0$  for   $x \in X^1_o$.
 Suppose that there exists a connected Euclidean open subset $U^1 \subset X^1_o$ (resp. $U^2 \subset X^2_o$) and an embedding $\varphi: U^1 \to U^2$ satisfying the following
two conditions. \begin{itemize} \item[(i)] $ \BP {\rm d} \varphi (\sC^1|_{U^1}) \subset \sC^2|_{\varphi(U^1)}$  and
\item[(ii)] for a general point $z \in U^1$ with $x = \varphi (z)$ and a general smooth point $\beta \in \sC^1_z$, the image  $\alpha:= \BP {\rm d} \varphi (\beta)$ is a smooth point of $\sC^2_x$ and
    $$\{ v \in T_{\alpha}(\sC^2_x), \ {\rm II}_{\sC_x^2, \alpha} (v, u) = 0 \mbox{ for all }  u \in {\rm d} (\BP {\rm d} \varphi) (T_{\beta}(\sC_z^1)) \} = 0.$$
 \end{itemize}
    Then there exists a rational map $\Phi: X^1 \dasharrow X^2$ such that $\Phi|_{U^1} = \varphi.$
 \end{theorem}
 To be precise, the statement in \cite{HongMok} is slightly weaker, although their argument essentially  proves   Theorem
 \ref{t.CFi}.  The full statement of Theorem \ref{t.CFi} is given in   \cite{Hwang14} with a simplified proof.

 The condition (ii) in Theorem \ref{t.CFi}  is a natural extension of the Gauss map condition in Theorem \ref{t.CartanFubini}. Note that the
 Cartanian condition in Theorem \ref{t.submanifold} on the characteristic connection is replaced by the algebro-geometric condition
 (ii).

 Theorem \ref{t.submanifold} goes beyond Theorem \ref{t.CartanFubini} even when $\dim X^1= \dim X^2$ because $\dim \sK^1$ can be strictly smaller than $\dim \sK^2$. For example,
  it can describe a finite morphism $X^1 \to X^2$ which sends members of $\sK^1$ to members of $\sK^2$.

  One immediate question is
 whether the map $\Phi$ can  be extended to a morphism $\Phi: X^1 \to X^2$ when $b_2(X^2) =1$, as in Theorem \ref{t.CartanFubini}.  In \cite{HongMok}, an affirmative answer is given when $X^1$ is modeled on a special class
 of Schubert submanifolds in rational homogeneous spaces $X^2=G/P$.

At present, submanifold theory in the interaction of Mori geometry and Cartan geometry is in an incipient stage.  Many natural questions can be raised, but the most fundamental one is an analogue of Problem \ref{p.central} in the submanifold setting. This would lead to a deeper aspect of Cartan geometry. Interesting applications to algebraic geometry are yet to come.

 \medskip
 In conclusion, VMRT-structures on uniruled projective manifolds provide us with a large number of examples in great diversity of cone structures admitting $\sP$-splitting characteristic connections. Cartan geometry of these structures has been understood for only a few special cases and a wide range of examples and problems remain to be explored.  The small number of results we have seen so far have already found interesting applications in algebraic geometry. Further development will undoubtedly bring more exciting applications. This will be a fertile ground for interactions between differential geometry and algebraic geometry.

\bigskip
{\bf Acknowledgment}.  I am very grateful to Ngaiming Mok and Richard Weiss for  valuable comments and helpful suggestions.



\end{document}